\documentclass[sn-mathphys,Numbered]{sn-jnl}

\usepackage{graphicx}%
\usepackage{multirow}%
\usepackage{amsmath,amssymb,amsfonts}%
\usepackage{amsthm}%
\usepackage{mathrsfs}%
\usepackage[title]{appendix}%
\usepackage{xcolor}%
\usepackage{textcomp}%
\usepackage{manyfoot}%
\usepackage{booktabs}%
\usepackage{algorithm}%
\usepackage{algorithmicx}%
\usepackage{algpseudocode}%
\usepackage{listings}%

\usepackage{todonotes}
\usepackage{enumitem}
\usepackage{mathtools}

\theoremstyle{thmstyleone}%
\newtheorem{example}{Example}[section]
\newtheorem{remark}{Remark}[section]
\newtheorem{definition}{Definition}[section]

\begin{document}

\title[Article Title]{Pragmatic Nonsense}


\author*[1]{\fnm{Ricardo} \sur{P. Cavassane}}\email{ricardo.peraca@gmail.com}

\author[1]{\fnm{Itala} \sur{M. Loffredo D’Ottaviano}}\email{itala@unicamp.br}
\equalcont{These authors contributed equally to this work.}

\author[1,2,3]{\fnm{Felipe} \sur{S. Abrahão}}\email{felipesabrahao@gmail.com}
\equalcont{These authors contributed equally to this work.}

\affil*[1]{\orgdiv{Centre for Logic, Epistemology and the History of Science}, \orgname{University of Campinas}, Brazil}

\affil[2]{ Oxford Immune Algorithmics, \orgname{Oxford University Innovation}, UK }

\affil[3]{ DEXL, \orgname{National Laboratory for Scientific Computing}, Brazil }



\abstract{Inspired by the early Wittgenstein’s concept of nonsense (meaning that which lies beyond the limits of language), we investigate two different types of nonsense: formal nonsense and pragmatic nonsense. The simpler notion of formal nonsense is defined in accordance with Tarski’s semantic theory of truth; the notion of pragmatic nonsense is in turn formulated within the context of the theory of pragmatic truth, also known as quasi-truth, as formalized by da Costa and his collaborators. Pragmatic nonsense extends formal nonsense, the same way da Costa’s pragmatic truth is an extension of Tarski’s definition of truth. An expression is thus considered formally nonsensical in case the formal criteria required for the assignment of any truth-value (whether true, false, pragmatically true, or pragmatically false) are not met; and an expression, or even a well-formed formula, is considered pragmatically nonsensical if either the formal or the pragmatic criteria of relevance (inscribed within the context of scientific practice) required for the assignment of any pragmatic truth-value (pragmatically true or pragmatically false) are not met. We also introduce the concept of strictly pragmatic truth, which excludes pragmatic nonsense and necessarily depends on certain criteria of relevance, unlike the original definition of pragmatic truth/quasi-truth.}




\maketitle

\section{Introduction}\label{sec1}

The distinction between sense and nonsense should be trivial in the context of formal languages, as will be clear in our definition of formal nonsense, initially formulated within the context of the semantic theory of truth of Alfred Tarski. However, as a formal language is equipped with formal tools to better represent the reality of a certain linguistic practice, such distinction may become more blurred. This occurs in the context of the theory of pragmatic truth of Newton da Costa and his collaborators, which aims at providing a formal account for the conception of truth assumed in the scientific practice. In order to clarify the distinction between sense and nonsense in that context, we define the notion of pragmatic nonsense. Notice that the pragmatic type of nonsense extends the formal type.

The early Wittgenstein’s concept of nonsense appears in the \textit{Tractatus Logico-Philosophicus}, originally published in German in 1922. We will follow its second English translation, by David Pears and Brian McGuinness \cite{bib1}.

The notion of formal nonsense will be initially defined within the context of Tarski’s semantic theory of truth, introduced in a 1933 paper originally published in Polish, the English translation with the title \textit{The concept of truth in formalized languages} published as a chapter of \textit{Logics, semantics, metamathematics} \cite{bib2}, and clarified in \textit{The semantic conception of truth and the foundations of semantics} \cite{bib3}. However, the definitions and nomenclature we will use here follow Mendelson’s \textit{Introduction to mathematical logic} \cite{bib4} and Hrbacek’s and Jech’s \textit{Introduction to set theory} \cite{bib5}. For the purposes of clarification, we will eventually contrast the terminology we use from the terminology found in other authors.

The notion of pragmatic nonsense, by its turn, is inscribed within the theory of pragmatic truth. The semantic-theoretical concept of pragmatic truth first appeared in \textit{Pragmatic truth and approximation to truth} \cite{bib6}, by Mikenberg, da Costa, and Chuaqui. That same year, it appears in da Costa’s \textit{Pragmatic Probability} \cite{bib7} with the name of quasi-truth, with a formalization later developed on \textit{The Logic of Pragmatic Truth} \cite{bib8}, which will be the one we will follow.

With our definition of pragmatic nonsense, we intend to clarify the pragmatic criteria necessary for the theory of pragmatic truth to achieve its goal, which is formalizing an intuitive notion of truth tacitly assumed in the scientific practice:

\noindent
\begin{quote}
[...] to accept a theory is to be committed, not to believing it to be true per se, but to holding it as if it were true, for the purposes of further elaboration, development and investigation. Thus acceptance involves belief that the theory is partially or pragmatically true only and this, we believe, corresponds to the fallibilistic attitude of scientists themselves \cite[p.~617]{bib8}.
\end{quote}

Therefore, the notion of formal nonsense will apply to all expressions that do not meet the formal criteria required for the assignment of any truth-value (whether true, false, pragmatically true, or pragmatically false) to them, mostly because they are not (well-formed) formulas; and the notion of pragmatic nonsense will apply also to all (well-formed) formulas that do not meet the pragmatic criteria for the assignment of any pragmatic truth-value (pragmatically true or pragmatically false) to them, pragmatic criteria which are inscribed within the context of scientific practice.

Thus, in Section ~\ref{sec2}, \textit{Truth and formal nonsense}, we will define the notion of formal nonsense for formal languages interpreted on total structures; then, in Section ~\ref{sec3}, \textit{Pragmatic truth and formal nonsense}, we will define the notion of formal nonsense for formal languages interpreted on simple pragmatic structures; and finally, in Section ~\ref{sec4}, \textit{Pragmatic nonsense}, we will the define the notion of pragmatic nonsense for formal languages interpreted on simple pragmatic structures.

\section{Truth and formal nonsense}\label{sec2}

In order to define nonsense in the context of formal languages, we must first assume a definition of formal language or, more specifically, of a first-order language (see \cite[p.~48-49]{bib4}), a definition of sentence of a first-order language (see \cite[p.~50]{bib4}), and a definition of interpretation of a first-order language (see \cite[p.~49]{bib4}). Since Mendelson \cite{bib4} defines an interpretation as consisting of a domain and assignments of relations to predicate symbols, operations to function symbols, and elements to individual constants, we must assume definitions of relation and operation or, in this case, function (see \cite[p.~19-29]{bib5}).

Assuming a definition of structure as a domain and a family of relations, we may say that a formal language is interpreted on a structure. Notice that Tarski and Vaught \cite{bib9} refer to the concept of structure with the term “relational system” or just “system”; in Hodges \cite{bib10}, by its turn, what we simply call structure is called “relational structure”. Notice also that what Mendelson \cite{bib4} calls “interpretation” has the name of “structure” in Enderton \cite{bib11}, in which the term “interpretation” has an entirely different meaning.

Then, we assume a definition of satisfaction (see \cite[p.~51-52]{bib4}) and a definition of truth of a formula in a first-order language (see \cite[p.~52]{bib4}).

None of the definitions listed above will be exposed here for they can be easily found, not only in the literature we refer to, but on many other basic sources.

Once a formal language is properly built and interpreted on a structure, it is possible to determine which sequences of symbols of such language, or expressions, have a determinate meaning (the closed formulas, or sentences, which are necessarily true or false), which expressions have an indeterminate meaning (the open formulas, which may come to be true or false), and which expressions have no meaning at all (for they are not considered formulas). A meaningless expression, therefore excluded from the well-formed formulas of a formal language, may be called nonsensical, if we use the term in the general sense that the early Wittgenstein gave it on the \textit{Tractatus Logico-Philosophicus}, as that which lies beyond the limits of language and has no possible truth-value \cite{bib1}; though we will not use the concept with the exact meaning it has in the \textit{Tractatus} (otherwise, this whole paper would be considered nonsense).

An example of a nonsensical expression would be “Socrates is identical”, or, in a first-order formal language with equality, `\textit{a} ='. ``The reason why ‘Socrates is identical’ means nothing is that there is no property called ‘identical'' \cite[~§5.473]{bib1}. Of course, Wittgenstein adds, “The proposition is nonsensical because we have failed to make an arbitrary determination, and not because the symbol, in itself, would be illegitimate” \cite[~§5.473]{bib1}, that is, when building our formal language, we defined the symbol ‘=’ as a substitute for the symbol `$ P_1^2 $', meaning thus a two-place (or dyadic) predicate of the form `$ P_1^2 (a, a) $' or `$ a = a $', not a one-place predicate, or property.

Following Wittgenstein’s conception of nonsense, in the context of Tarski’s semantic theory of truth, the notion of formal nonsense would be simply defined as an expression which fails to be a (well-formed) formula.

\begin{definition}[\textbf{Formal nonsense for formal languages interpreted on total structures}]
\label{def:21}
Let $ L $ be a formal language. An expression `$ \alpha $' of $ L $ is \emph{formally nonsensical} if, and only if, ‘$\alpha$’ is not a (well-formed) formula.
That is, given an interpretation $ \mathbf{I} $ on a structure $ E $ with a domain $ D $ according to which $ L $ is interpreted, ‘$\alpha$’ cannot be interpreted.
\end{definition}

This definition is rather trivial, and will remain quite simple when we extend it in the next section, after defining pragmatic truth (see Definition \ref{def:36}). However, when we enter the domain of the theory of pragmatic truth, it will become evident that a stronger notion of nonsense, which we will call pragmatic nonsense (see Section~\ref{sec4}), will be necessary in order to prevent a weakening of pragmatic truth, therefore necessarily giving rise to a stricter notion of pragmatic truth in which subjective or informal criteria of relevance is taken into account.

\section{Pragmatic truth and formal nonsense}\label{sec3}

As is the case with the Tarskian definition of true sentence, the definition of pragmatically true sentence only applies to sentences of formalized languages. Thus, the definition of formal language assumed on the previous section will be maintained, as well as every other previously assumed definition, with the exception of the definitions of relation and structure, which must be replaced with the definitions of partial relation and partial structure, with no loss, since a total relation is a particular case of partial relation and a total structure is a particular case of partial structure. 
For the sake of clarifying the forthcoming concepts and definitions, we briefly restate the definitions of partial relation and partial structure. 

\begin{definition}[\textbf{Partial relation}]
Let $ D $ be a non-empty set. An \textit{$ n $-ary partial relation $ R^n $} on $ D $ is an ordered triple $ \left< R_1, R_2, R_3 \right> $, in which $ R_i \cap R_j = \emptyset $, for $ i  \neq j $, $ i $, $ j \in \left\{ 1, 2, 3 \right\} $, and $ R_1 \cup R_2 \cup R_3 = D_n $, such that:
\begin{enumerate}
\item $ R_1 $ is the set of the ordered \textit{n}-tuples of elements of $ D $ that we know that belong to $ R^n $;

\item  $ R_2 $ is the set of the ordered \textit{n}-tuples of elements of $ D $ that we know that do not belong to $ R^n $;

\item $ R_3 $ is the set of the ordered $ n $-tuples of elements of $ D $ that we do not know if they belong to $ R^n $ or not.
\end{enumerate}
\end{definition}

\begin{definition}[\textbf{Partial structure}]
A \emph{partial structure} $ E $ is an ordered pair $ \left< D, R_i^j \right>_{ i \, \in \, I, \, j \, \in \, J }$, such that:

\begin{enumerate}
\item $ D $ is a non-empty set, the domain of the structure;

\item $ \left( R_i^j \right)_{ i \, \in \, I, \, j \, \in \, J }$, such that $ J = \left\{ 1 , 2 , \dots , \textit{n} \right\} $ and $ I $ is a family of $ j $-ary partial relations $ R_i^j $ on $ D $.
\end{enumerate}

\end{definition}

\begin{remark}[Total relation and total structure]
If, for a given partial relation $ R^n $, $ R_3 = \emptyset $ holds, then $ R^n $ is a total $ n $-ary relation, or just a relation in the usual sense. 

(Notice that a different definition of partial structure may be found in \cite{bib6}).

If, for every partial relation of a given partial structure $ E $, $ R_3 = \emptyset$, then $ E $ is a total structure, or just a structure in the usual sense.
\end{remark}

Having defined the notions of partial relation and partial structure, in order to define pragmatic truth we must then restate the definitions of simple pragmatic structure (SPS) and $ \mathcal{A} $-normal structure, as well as the definition of expansion, originally defined in \cite{Bueno1996} for partial models, but here defined for partial structures. 


\begin{definition}[\textbf{Expansion of a structure}]
Let $ E_1 $ be a partial structure and $ E_2 $ a partial or total structure.
$ E_2 $ \textit{expands} (or is an \textit{expansion} of) $ E_1 $ if, and only if:

\begin{enumerate}

\item The domain $ D $ of $ E_1 $ is the same domain of $ E_2 $;


\item For each partial relation $ R_i^j = \left< R_1 , R_2 , R_3 \right> $ of $ E_1 $ there is a corresponding partial or total relation $ R_k^j  = \left< R'_1, R'_2 , R'_3  \right> $ of  $ E_2 $ such as that:

\begin{enumerate}

\item $ R_1 \subseteq R'_1 $;

\item $ R_2 \subseteq R'_2 $;

\item $ R_3 \subseteq R'_1  \cup R'_2 \cup R'_3 $.

\end{enumerate}

\end{enumerate}

\end{definition}


\begin{definition}[\textbf{Simple pragmatic structure}]
A \textit{simple pragmatic structure} (SPS) $ \mathcal{A} $ is an ordered triple $ \left< D , R_i^j , \mathcal{P} \right>_{ i \, \in \, I, \, j \, \in \, J } $, such that:
\begin{enumerate}
\item $ D $ is a non-empty set, the domain of the structure;

\item $ \left( R_i^j \right)_{ i \, \in \, I, \, j \, \in \, J } $, such that $ J = \{ 1, 2, \dots , n \} $ and $ I $ is a family of $ j $-ary partial relations $ R_i^j $ on $ D $;

\item $ \mathcal{P} $ is a possibly empty set of sentences of a formal language $ L $ interpreted by the interpretation $ \mathbf{I} $ on $ \mathcal{A} $, the set of primary sentences of $\mathcal{A}$.
\end{enumerate}

\end{definition}

\begin{definition}[\textbf{$\mathcal{A}$-normal structure}]
Let $ \mathcal{A} $ be a simple pragmatic structure such that $\mathcal{A} = \left< D, R_i^j, \mathcal{P} \right>_{ i \, \in \, I, \, j \, \in \, J } $, $\mathcal{B} $ a total structure such that $ \mathcal{B} = \left< D , R_k^j, \mathcal{P} \right>_{ k \, \in \, K, \, j \, \in \, J } $, and $ L $ a formal language interpreted by the interpretation $ \mathbf{I} $ on $\mathcal{B} $ and on $ \mathcal{A} $. 
$ \mathcal{B} $ is an $ \mathcal{A} $-\emph{normal structure} if, and only if:
\begin{enumerate}

\item $ \mathcal{B} $ is an expansion of $ \mathcal{A} $;

\item  If `$ \alpha $' is a sentence of $ L $ that belongs to $ \mathcal{P} $, then `$ \alpha $' is true according to the interpretation $ \mathbf{I} $ of $ L $ in $ \mathcal{B} $.
\end{enumerate}

\end{definition}

\begin{remark}[The existence of $ \mathcal{A} $-normal structures]
In general terms, the existence of $ \mathcal{A} $-normal structures is tied to the existence of a model for the consequence set $ Cn\left( \mathcal{P} \right) $ obtained from all possible derivations of $ \mathcal{P} $ together with the set of all closed literals mapping each tuple in the partial relations $ R_i^j $ into the language of the $ \mathcal{A} $.
The necessary and sufficient conditions for the existence of at least one $ \mathcal{A} $-normal structure, given a SPS, can be found in \cite{bib6}.
\end{remark}

\begin{remark}[Partial models]
Notice that a SPS is usually defined for relational structures, but it can be straightforwardly understood as a model by taking function symbols as a particular case of relations/predicates for which there is only a single rightmost element for any other elements in the tuple, and the constant symbols are just unary relations/predicates.
\end{remark}

With the above definitions, one is now able to set the conditions in which a sentence becomes pragmatically true and those in which a sentence becomes pragmatically false. Notice that, after defining pragmatic nonsense in Section~\ref{sec4}, we will refine the definition of pragmatic truth below
with the definition of \emph{strictly} pragmatic truth. 

\begin{definition}[\textbf{Pragmatic truth and pragmatic falsehood of formulas}]
Let $ \mathcal{A} $ be a simple pragmatic structure, $\mathcal{B} $ an $ \mathcal{A} $-normal structure, $ L $ a formal language interpreted by the interpretation $ \mathbf{I} $ on $ \mathcal{A} $ and on $ \mathcal{B} $, and `$ \alpha $' a formula of $ L $.
\begin{itemize}
\item A formula `$ \alpha $' is \emph{pragmatically true} (or \emph{quasi-true}) according to the interpretation $ \mathbf{I} $ in a simple pragmatic structure $\mathcal{A} $ if, and only if, `$ \alpha $' is true according to the interpretation $ \mathbf{I} $ in an $\mathcal{A} $-normal structure $ \mathcal{ B }_i $ (i.e., if, and only if, there exists at least one $ \mathcal{A} $-normal structure $ \mathcal{ B }_i $ in which `$ \alpha $' is true).

\item A formula `$ \alpha $' is \emph{pragmatically false} (or \textit{quasi-false}) according to the interpretation $ \mathbf{I} $ in a simple pragmatic structure $ \mathcal{A} $ if, and only if, `$ \alpha $' is not pragmatically true according to the interpretation $ \mathbf{I} $ in $ \mathcal{A} $ (i.e., if, and only if, there does \emph{not} exist even one $\mathcal{A} $-normal structure $ \mathcal{B}_i $ in which `$ \alpha $' is true).

\end{itemize}
\end{definition}

Let us see an example, in which a family of four (a father, a mother, and two sons) and some of their relationships are described as a partial structure, and then as a simple pragmatic structure and an $ \mathcal{A} $-normal structure, which will allow us to state pragmatic truths about that family.

\begin{example}[Partial structure of a family $ F $]
\label{ex:1}
The partial structure $ F $ is an ordered pair $ \left< D_F , R_i^2 \right>_{ i \, = \, \left\{ 1 , 2 \right\} } $, with $ R_i^2 = \left\{ R_1^2 , R_2^2 \right\} $, such that:

\begin{enumerate}
\item $ D_F = \{ $Joseph, Mary, John, Peter$ \} $;

\item $ R_1^2 = $ ⟨\textit{$R_1$, $R_2$, $R_3$}⟩, such that: $R_1 = $ {⟨Joseph, Mary⟩, ⟨Mary, Joseph⟩}; $R_2 = $ {⟨Joseph, Joseph⟩, ⟨Joseph, John⟩, ⟨Joseph, Peter⟩, ⟨Mary, Mary⟩, ⟨Mary, John⟩, ⟨Mary, Peter⟩, ⟨John, Joseph⟩, ⟨John, Mary⟩, ⟨John, John⟩, ⟨John, Peter⟩, ⟨Peter, Joseph⟩, ⟨Peter, Mary⟩, ⟨Peter, John⟩, ⟨Peter, Peter⟩}; and $ R_3 = \emptyset $;

\item $ R_2^2 = $ ⟨\textit{$R_1$, $R_2$, $R_3$}⟩, such that: $ R_1 = $ {⟨Joseph, John⟩}; $R_2 = $ {⟨Joseph, Joseph⟩, ⟨Joseph, Mary⟩, ⟨Mary, Joseph⟩, ⟨Mary, Mary⟩, ⟨Mary, John⟩, ⟨Mary, Peter⟩, ⟨John, Joseph⟩, ⟨John, Mary⟩, ⟨John, John⟩, ⟨John, Peter⟩, ⟨Peter, Joseph⟩, ⟨Peter, Mary⟩, ⟨Peter, John⟩, ⟨Peter, Peter⟩}; and $ R_3 = $ {⟨Joseph, Peter⟩}.
\end{enumerate}

\end{example}

So, we know that Joseph is married to Mary, and that is represented by the relation $ R_1^2 $ (which, since $ R_3 $ is empty, is a total relation); and we also know that Joseph is the biological father of John, but we do not know if Joseph is the biological father of Peter or not, and that is represented by the relation $ R_2^2 $ (which is, thus, a partial relation).

Consider, in the following example, built upon Example \ref{ex:1}, that a DNA paternity test was made in order to answer that question, and that both the result of the test and the implication of such result were added to the set of primary sentences of a simple pragmatic structure.

\begin{example}[Simple pragmatic structure $ \mathcal{G} $]
\label{ex:2}
The simple pragmatic structure $ \mathcal{G} $ is an ordered triple $ \left< D_\mathcal{G} , R_i^j , \mathcal{P} \right>_{ i \, \in \, I, \, j \, \in \, J }$, such that:

\begin{enumerate}
\item $ D_\mathcal{G} $ is a non-empty set, the domain of the  structure; 

\item $ \left( R_i^j \right)_{ i \, \in \, I, \, j \, \in \, J } $, such that $ J = \left\{ 1 , 2 , \dots , \textit{n} \right\} $ and $ I $ is a family of $ j $-ary partial relations $ R_i^j $ on $ D_\mathcal{G} $; 

\item $ \mathcal{P}_\mathcal{G} $ is a set of sentences of a formal language $ L_\mathcal{G} $ interpreted by the interpretation $ \mathbf{I}_{ \mathcal{G} } $ on $\mathcal{G} $, such that $ \mathcal{P}_\mathcal{G} = \left\{ \alpha, \alpha \to \beta \right\} $ (`$\alpha$' may be translated, in the metalanguage $ M_\mathcal{G} $ of $ L_\mathcal{G} $, as `Joseph has at least $99,99$\% of chance of being the biological father of Peter'; `$\beta$' may be translated, in the metalanguage $ M_\mathcal{G} $ of $ L_\mathcal{G} $, as `Joseph is the biological father of Peter'; and, therefore, `$ \alpha \to \beta $' may be translated, in the metalanguage $ M_{ \mathcal{G} } $ of $ L_\mathcal{G} $, as `If Joseph has at least $99,99$\% of chance of being the biological father of Peter, then Joseph is the biological father of Peter').
\end{enumerate}

\end{example}

Now, consider that the simple pragmatic structure in question was expanded to a total structure, forming an $\mathcal{A}$-normal structure.

\begin{example}[$ \mathcal{A} $-normal structure $ \mathcal{H} $ for Example~\ref{ex:2}]
The $ \mathcal{A} $-normal structure $ \mathcal{H} $ is a total structure such that:
\begin{enumerate}

\item $ \mathcal{H} $ is an expansion of $ \mathcal{G} $ (the partial relation $ R_2^2 = \left< R_1 , R_2 , R_3 \right> $ of $ \mathcal{G} $ has a correspondent total relation $ { R' }_2^2  = \left< R'_1 , R'_2 , R'_3 \right> $ in $ \mathcal{H} $ such that $ R'_1  = $ {⟨Joseph, John⟩, ⟨Joseph, Peter⟩});

\item  The sentences `$ \alpha $' and `$ \alpha  \to \beta $' of $ L_{ \mathcal{G} } $ belonging to $\mathcal{P}_{ \mathcal{G} } $ are true according to the interpretation $ \mathbf{I}_{\mathcal{G}}$ of $ L_{\mathcal{G}} $ in $ \mathcal{H} $.
\end{enumerate}

\end{example}

Equipped with the partial structure $ F $, the simple pragmatic structure $ \mathcal{G} $, and the $\mathcal{A} $-normal structure $ \mathcal{H} $ (and, of course, with a formal language $ L_{\mathcal{G}} $ interpreted by an interpretation $ \mathbf{I}_{ \mathcal{G} } $ on $ F $, $ \mathcal{G} $ and $\mathcal{H} $), one can assert, for instance, that the sentence `$ P_2^2 $ (\textit{Joseph}, \textit{Peter})' or, in a metalanguage $ M_{ \mathcal{G} } $, `Joseph is the biological father of Peter', is pragmatically true in the simple pragmatic structure $ \mathcal{G} $, since it is true in the $ \mathcal{A} $-normal structure $\mathcal{H} $. That is possible even though one does not know if that sentence is true or false in the partial structure $ F $, since the sentence is consistent with the primary sentences of $ \mathcal{P}_{ \mathcal{G} } $. In this example, it is not only consistent with but also a logical consequence of the sentences of $\mathcal{P}_{ \mathcal{G} } $, but that is due to the simplicity of the example and must not always be the case; also, in this example, only one $ \mathcal{A} $-normal structure is possible out of the two possible expansions of the simple pragmatic structure, but that is also due to the simplicity of the example.

\subsection{The set of primary sentences}\label{sectionPrimarySentences}

Before we proceed to the definition of pragmatic nonsense, we must first elucidate the role of the set of primary sentences $ \mathcal{P} $ of a simple pragmatic structure. The sentences of $ \mathcal{P} $ constitute an empirical and/or theoretical framework which provides a foundation for pragmatic truth. Such primary sentences ``[...] are based on observation and experience [...] or were instituted by previous investigations [...]'' \cite[p.~204]{bib6}) and ``[...] in general include both [...] true decidable sentences, such as observation statements, and certain general propositions encompassing laws already assumed to be true'' \cite[p.~162]{bib12}. They may also constitute abductive hypotheses, and may be either true in the Tarskian sense, pragmatically true themselves, or ``accepted'' as true. 

The sentences of $ \mathcal{P} $ have, at the same time, the function of limiting the possible $ \mathcal{A} $-normal structures a simple pragmatic structure may be expanded to, for those expansions that are inconsistent with the sentences of $ \mathcal{P} $ are not considered $ \mathcal{A} $-normal and thus disregarded.
This highlights both the motivation and applicability of the concept of partial structures and quasi-truth to model theory \cite{Bueno1996}.
In addition to the notion of those statements beyond what one can achieve by inferences alone, quasi-truth brings this \emph{partialness} and \emph{openness} to semantics, i.e., to whether something has a model or is satisfiable.
It instantiates the gain of irreducibly new truths that occurs for deducibility, but also for satiafiability: 
as one adds a novel statement, which was previously unprovable by the current theory, the unknown or undetermined components in $ R_3 $ are also reduced along with the reduction in the number of unprovable sentences.
This is shown in the following two Examples~\ref{exampleClearcut} and~\ref{exampleGradual}.

\begin{example}\label{exampleClearcut}
Suppose the traditional set $ \mathcal{S} $ of axioms of Peano Arithmetics (PA) is consistent, $ L_{ \mathcal{G} } $ is the first-order language of PA, and $ \mathcal{G} = \left< \textit{D}_{\mathcal{G}}, R_i^j, \mathcal{P} \right> $ is a simple pragmatic structure defined upon the set $ \mathcal{P} = \mathcal{S} $ of primary sentences.
This structure $ \mathcal{G} $ has the same\footnote{Thus, some of the relation/predicate symbols are in fact restricted to the case of function symbols or constant symbols.} signature and atomic diagram \cite{bib10} of the standard model $ \mathfrak{N} $ of PA, except for an additional dyadic \emph{partial} relation $ R_{ Pf } $, which is added to the signature of $ \mathfrak{N} $ in order to construct the signature of $ \mathcal{G} $, so that $ \left( x , y \right) \in R_{ Pf } $ means that `$ x $ is the Gödel number in $ \mathfrak{N} $ of a proof of the sentence encoded by the number $ y $ in $ \mathfrak{N} $'.

Notice that one shall not confuse the (partial) relation $ R_{ Pf } $ with the notation $ Pf_\mathcal{S}(x,y) $ (or $ Prov_\mathcal{S}\left( x , y \right) $) usually seen in the literature to abbreviate the well-formed formula in $ L_{ \mathcal{G} } $ that represents the (recursive) relation `$ x $ is the Gödel number of a proof of the sentence encoded by the number $ y $' \cite{bib4,bib11,Smorynski1977}.
As usual, let $ \text{CON}({\mathcal{S}}) $ denote the sentence that states the consistency of $ \mathcal{S} $ such as the formula $ \text{CON}(S) \coloneqq \text{` } \neg \exists x Pf_\mathcal{S}\left( x , 0 \neq 0 \right) $' \cite{Smorynski1977}.

Regarding the atomic partial\footnote{ Being composed only of partial relations or functions, an atomic partial diagram is defined by those atomic sentences that are mapped into the first components $ R_1 $'s of each partial relation/function in the structure.} diagram of $ R_{ Pf } $, as previously stated for the meaning of $ \left( x , y \right) \in R_{ Pf } $:
\begin{enumerate}
\item the first component $ R_1 $ of the partial relation $ R_{ Pf } = \left< R_1 , R_2 , R_3 \right> $ is defined upon all, and only, atomic sentences true in the standard model (i.e., $ R_1 $ is composed of all pairs of natural numbers whose first number encodes a correct proof in PA and the second natural number encodes a sentence in $ L_{ \mathcal{G} } $); 
\item however, unlike the traditional total-relation case, its second component $ R_2 $ is defined upon all sentences that are false in every model of PA, together with those that are false in at least one model of $ \text{PA} \, + \, \neg \text{CON}({\mathcal{S}}) $ but also true in at least another nonstandard model of PA;
\item and its third component $ R_3 $ is defined upon all, and only, the atomic sentences that are true in every model of $ \text{PA} \, + \, \neg \text{CON}({\mathcal{S}}) $.
\end{enumerate}
Thus, notice that, except for the relation $ R_{ Pf } $, $\mathcal{G} $ is isomorphic to the standard model $ \mathfrak{N} $ of PA.
One can understand the partial structure $ \mathcal{G} $ as the one capturing the intuition of the expected model for what is provable in PA.
One already has formal knowledge (instantiated in the first component $ R_1 $) sufficient for accounting for everything already provable in the standard model $ \mathfrak{N} $.
For those pairs of numbers upon which every model agrees about the non-consistency of PA, there is a lack of formal knowledge (instantiated in the third component $ R_3 $), since in this case they apparently should be considered part of any model (due to the completeness of the first-order logic), but ultimately fail at least for one model (i.e., the standard model).

Now, let 
\begin{equation*}
	\mathcal{P}' = \mathcal{P} \cup \left\{ \text{CON}\left( \mathcal{S} \right) , \, \forall x , y \left( Pf_\mathcal{S}\left( x , y \right) \leftrightarrow \text{`$ R_{ Pf } $'}\left( x , y \right) \, \right) \right\}
	\text{ ,}
\end{equation*}
where $ L_{ \mathcal{G}' } $ is the new language of PA augmented\footnote{ Possibly, if necessary, it can be also augmented by a denumerable number of additional constant symbols.} by the new relation symbol $ \text{`$ R_{ Pf } $'} $. 
That is, in order to build $ \mathcal{P}' $, one adds the sentence ‘$ \text{CON}({\mathcal{S}}) $’ to the set $ \mathcal{P} $ together with the sentence stating the equivalence between the encoding of a proof in arithmetics and the new relation symbol $ \text{`$ R_{ Pf } $'}\left( x , y \right) $.
Therefore, one can directly demonstrate that there will be \emph{no} total ($ \mathcal{G}' $-normal) structure expanding $ \mathcal{G}' = \left< \textit{D}_{\mathcal{G}}, R_i^j, \mathcal{P}' \right> $ in which a pair from $ R_3 $ is added to $ R_1 $.
In this case, every atomic sentence $ \text{`$ R_{ Pf } $'}\left( n , m \right) $ corresponding to a pair in $ R_3 $ will need to then correspond to a pair added to $ R_2 $.
This is because $ R_3 $ contains those pairs that are part of any model of $ Cn\left( \mathcal{P} \cup \left\{ \exists x Pf_\mathcal{S}\left( x , 0 \neq 0 \right) \right\} \right) $, pairs which were not already included in $ R_1 $ by construction.
As a consequence, adding those pairs to $ R_1 $ would automatically prevent $ \mathcal{G}' $ from having a total structure as a model that satisfies $ Cn\left( \mathcal{P}' \right) $.
The reader is invited to note that instead of $ \mathcal{S} $, the argument presented in this Example~\ref{exampleClearcut} applies analogously to every formal theory $ \mathbf{T} $ strong enough to contain Peano Arithmetics.
\end{example}

The addition of a couple of sentences to $ \mathcal{P} $ in the above Example~\ref{exampleClearcut} gives rise to a semantic ``clear cut'' into a partial structure (in the case, $ \mathcal{G} $) undergone by the total structures (that expands $ \mathcal{G} $), so that every unknown pair collapses (into falsehood, i.e., $ R_2 $), and thus the unknowability disappears completely.
The partial relation $ R_{ Pf } $ was defined with the purpose of making the consistency of PA trigger such a collapse.
We posit that the existence of such ``clear-cut'' formal theories constitutes a fruitful line of future research.
The problem is to investigate the conditions in which simple pragmatic structures would necessarily collapse into a total structure if one adds a \emph{finite} number of sentences to the theories those structures are meant to satisfy.

As shown in the following Example~\ref{exampleGradual}, pragmatic truth can also reflect a \emph{gradual} loss of unknowability as one adds new unprovable sentences to $ \mathcal{P} $.
This gradual gain of knowledge in the partial structures that model a set of sentences can follow along the very gain in deducibility power.
More formally, the partial structures for the formal theories $   \mathbf{T} $ (i.e., the set $ \mathcal{P} $ of primary sentences) in Example~\ref{exampleGradual} are defined in such a way that they are restricted to $ Cn\left( \mathbf{T} \right) $, i.e., to what can be deduced from $ \mathcal{P} $.
In consonance with the purpose of being a formalization of a pragmatist theory of truth that encompasses both the correspondence and the coherence proposals \cite{bib12,Bueno1996,DOttaviano2007}, Example~\ref{exampleGradual} highlights how the unknowability present in the syntactical realm can be reflected in the semantical realm, and vice versa.

\begin{example}\label{exampleGradual}
Let $ \mathcal{S} $ be the traditional set of axioms of Peano Arithmetics (PA) on which we suppose it is consistent. 
Let $ L_{ \mathcal{G} } $ be the first-order language of PA, and $ \mathcal{G}_k = \left< \textit{D}_{\mathcal{G}}, R_i^j, \mathcal{P} \right> $, where $ k \geq 0 $, be a simple pragmatic structure defined upon the set $ \mathcal{P} = \mathcal{S} $ of primary sentences whose signature is exactly the same as the standard model $ \mathfrak{N} $ of PA.
Notice that as in the previous Example~\ref{exampleClearcut}, we are only using the predicate notation $ R_i^j $ for the sake of simplifying notation, although in fact the structure also includes function symbols and constant symbols.
Let $ \vdash_k $ with $ 0 \leq k < \infty $ denote the deducibility in at most $ k $ proof steps.

Let $ C \in \mathbb{N} $ be an arbitrarily large constant.
Unlike in the Example~\ref{exampleClearcut}, the atomic partial diagram of each $ R_i^j = \left< R_1 , R_2 , R_3  \right> $ (including the relations and functions `$ = $', `$ < $', `$ + $', `$ \times $', `$ Exp $', `$ S $') of the SPS $ \mathcal{G}_k $ is defined by:
\begin{enumerate}
	\item the first component $ R_1 $ is defined by at most $ k + C $ tuples that \emph{belong} to $ \mathfrak{N} $ such that:
	 \begin{enumerate}
		 \item they were already included in $ \mathcal{G}_m $ with $ m < k $;
		 
		 \item there is at least one sentence $ \sigma $ and at least one element of the tuple on which a variable $ x $ of the sentence $ \sigma $ is interpreted  such that
		$ \mathcal{P} \vdash_k \sigma $, where $ \sigma $ is true in $ \mathfrak{N} $, $ \sigma = \exists x \phi\left( x \right) $, and $ \phi\left( x \right) $ is a well-formed formula;


		  \item or iff $ 1.2 $ does not apply, $ \sigma $ is an axiom of PA, where there is at least one element of the tuple on which a variable $ x $ of the sentence $ \sigma $ is interpreted. 
	 \end{enumerate}

	\item the second component $ R_2 $ is defined by all, and only, tuples that do \emph{not} belong to any model of PA;
	
	\item the third component $ R_3 $ comprising all, and only, tuples that were not included in $ R_1 $ or $ R_2 $.
	
\end{enumerate}

Intuitively, $ \mathcal{G}_k $ is a partial structure/model for PA that mirrors the deducibility capabilities of $ \mathcal{S} $ into the structure itself, that is, it ``forces'' or ensures that semantics follow along syntactics.
As the value of $ k $ increases, more tuples can be included in the first components of the partial relations in addition to those already included for any $ \mathcal{G}_{ m } $ with $ m < k $.
If $ \zeta $ is an element of a tuple in a third component of a partial relation of $ \mathcal{G}_k $ for which the variable $ x $ is interpreted as $ \zeta $, then this tuple may eventually be included to $ R_1 $ for some $ \mathcal{G}_{ k' } $ with $ k' > k $.

We have that $ \mathcal{G}_{ k + 1 } $ would first add those new tuples containing at least one element for which a newly proven sentence $ \exists x \phi\left( x \right) $ appears.
Only if no such a sentence appears is that another tuple from the standard model is chosen in lexicographical order to be added to $ R_1 $ of the respective relation.

Clearly, $ \mathcal{G}_0 $ includes all the first $ 6 \left( 0 + C \right) $ tuples (in lexicographical order) that belong to the standard model of PA.
This is because we assumed that PA was consistent, hence satisfiable by some $ \mathfrak{N} $.
By construction in condition $ 1.3 $, it is guaranteed that even if no formula $ \phi\left( x \right) $ is ever found, $ 6 $ tuples of $ \mathfrak{N} $ are always added each time $ k $ increases by $ 1 $ (e.g., $ x $ being one of the variables occurring in the axiom $ \forall x \left( S(x) \neq 0 \right) $). 
From the compactness of first-order logic \cite{bib11} and the computability of the atomic diagram \cite{Cenzer2020StructureRandomnessComputability} of $ \mathfrak{N} $, one can directly demonstrate that the limit structure of this process, denoted by $ \mathcal{G}_\omega $ (where $ \omega $ is the usual first infinite ordinal) will be a partial structure isomorphic to $ \mathfrak{N} $, except for the third components in each partial relation.
For example, notice that the element in the total relations of a nonstandard model of PA that would satisfy $ \neg \text{CON}(S) = \exists x Pf_\mathcal{S}\left( x , 0 \neq 0 \right) $ will never be added to the first or second components of the relations of $ \mathcal{G}_\omega $.
Thus, to finally convert $ \mathcal{G}_\omega $ into a total structure one might still need to move everything else left in each $ R_3 $ into $ R_2 $, respectively.

One can also demonstrate that the standard model $ \mathfrak{N} $ is one of the possible $ \mathcal{G}_k $-normal models of any $ \mathcal{G}_k $ or $ \mathcal{G}_\omega $.
However, if one adds $ \neg \text{CON}({\mathcal{S}}) $ to $ \mathcal{P} $, then at least one new tuple not yet included in $ \mathfrak{N} $ will need to be included into any  $ \mathcal{G}'_k $-normal model, where $ \mathcal{G}'_k = \left< \textit{D}_{\mathcal{G}}, R_i^j, \mathcal{P}' \right> $ and $ \mathcal{P}' = \mathcal{P} \cup \left\{ \neg \text{CON}\left( \mathcal{S} \right) \right\} $, therefore changing the possible total structures that expand the partial structure.
\end{example}

In general, even though we may not know exactly everything, whether it is or not the case, in a given partial structure, we may know what may or may not be the case given certain data previously obtained, certain theoretical results previously accepted, or simply some investigation hypotheses.
However, from the formal point of view, ${\mathcal{P}}$ may be empty; in principle, there are also no restrictions or impositions for the inclusion of sentences in ${\mathcal{P}}$, except that they must be sentences of a formal language $ L $ interpreted by an interpretation $ \mathbf{I} $ on the referred simple pragmatic structure. 
With the introduction of the notion of pragmatic nonsense, the present article aims at introducing pragmatic criteria capable of clarifying the inclusion of sentences in ${\mathcal{P}}$.

\subsection{Formal nonsense in simple pragmatic strutures}\label{sectionFormalNonsenseSPS}

Before we move on to the notion of pragmatic nonsense, though, we must update our definition of formal nonsense one last time. 
For even with the introduction of the tools of the theory of pragmatic truth, there remains a case in which a (well-formed) formula will be considered formally nonsensical: when a new sentence puts the existence of at least one totalizing expansion of the simple pragmatic structure in question, that is, adding a sentence to the set $ \mathcal{P} $ of assumedly true sentences would render the existence of an ${\mathcal{A}}$-normal structure not determined.

\begin{definition}[\textbf{Formal nonsense for formal languages interpreted on simple pragmatic structures}]
\label{def:36}
Let $ L $ be a formal language interpreted according to an interpretation $ \mathbf{I} $ with a domain $ D $, and a simple pragmatic structure $ \mathcal{A} = \left< D , R_i^j , \mathcal{P} \right> $. An expression `$ \alpha $' of $ L $ is \emph{formally nonsensical} if, and only if: either `$ \alpha $' is not a (well-formed) formula; or, in case `$ \alpha $' is a (well-formed) formula, one cannot determine the existence of at least one $ \mathcal{A}' $-normal structure $ \mathcal{B} $ built upon $ \mathcal{A}' = \left< D , R_i^j , \mathcal{P} \cup \left\{ \alpha \right\} \right> $.
\end{definition}

Notice that Definition~\ref{def:36} applies analogously for a scheme of formulas $ \overline{ \alpha }  $ instead of a single formula $ \alpha $.

Besides the usual formal nonsense as in Definition~\ref{def:21}, the one for SPSs in Definition~\ref{def:36} includes other less trivial examples of nonsense.
In consonance with what one might expect, one example of SPS formal nonsense would happen when one tries to force the predicate of truths back into the language.
Indeed, in Example~\ref{exampleClearcut} we saw that if PA is consistent, one can force the predicate symbol for provability into the language of PA \cite{Smorynski1977}.
However, for sufficiently strong languages and theories, this method would lead to contradictions, i.e., impasses that would render the addition of such formulas nonsense.
This is shown in Example~\ref{exampleFormalnonsense} below.

\begin{example}\label{exampleFormalnonsense}
Suppose a sufficiently strong and satisfiable theory $ \mathcal{T} $ that encompasses Peano Arithmetics (PA), any arbitrarily chosen set theory such as Zermelo–Fraenkel-Choice (ZFC), and every possible acceptable proof in model theory.
Let $ L $ be a sufficiently expressive computable language that includes every sentence in the language of $ \mathcal{T} $ plus the additional expression `$ \text{TRUTH}\left( \cdot \right) $'.
Let $ \mathcal{G} = \left< \textit{D}_{\mathcal{G}}, R_i^j, \mathcal{P} \right> $ be a (total) pragmatic structure defined upon the set $ \mathcal{P} = \mathcal{T} $ of primary sentences such that this structure is a model for $ \mathcal{T} $ and $ Cn\left( \mathcal{T} \right) $ in the augmented language $ L $. 
Thus, since it is a total structure, any third component of any relation $ R_i^j $ in $ \mathcal{G} $ is empty.

Let $ \text{TRUTH} = \left< R_1 , R_2 , R_3 \right> $ be a unary partial relation/predicate such that:
\begin{enumerate}
	\item $ x \in R_1 $ iff $ x $ is the Gödel number of a sentence true in 
	$ \mathcal{G} $
	or in the pragmatic structure 
	\[ \left< \textit{D}_{\mathcal{G}}, R_i^j + \text{TRUTH} , \mathcal{P} \cup \left\{ \text{`TRUTH}\left( x \right) \text{'} \, \middle\vert \, x \in R_1  \right\} \cup \left\{ \neg \text{`TRUTH}\left( x \right) \text{'} \, \middle\vert \, x \notin R_1  \right\} \right> \text{ ;} \]
	
	\item $ x \in R_2 $ iff $ x $ is the Gödel number of a sentence that is false in any model of $ \mathcal{T} $;
	
	\item $ x \in R_3 $ iff $ x $ is the Gödel number of a sentence not yet included in $ R_1 $ or $ R_2 $.
\end{enumerate}

Because we are dealing with partial relations and not total relations, one can always move into $ R_3 $ every Gödel number of a sentence in $ L $ for which one cannot determine whether it should be included in $ R_1 $ or not---e.g., due to generating self-referential contradictions---, in this manner rendering $\text{TRUTH}$ a well-defined \emph{partial} relation, although never being a well-defined \emph{total} relation.

Notice that since no expression of $ \mathcal{T} $ or $ Cn\left( \mathcal{T} \right) $ is ever interpreted in $ \text{TRUTH} = \left< R_1 , R_2 , R_3 \right> $, then adding this (partial) relation to the structure $ \mathcal{G} $ would not change any satisfiability condition of any sentence in $ Cn\left( \mathcal{T} \right) $. 
In other words, if $ R_1 $ is well-defined, the SPS $ \mathcal{G}' = \left< \textit{D}_{\mathcal{G}} , R_i^j + \text{TRUTH} , \mathcal{P} \right> $ would also be a (partial) model of $ \mathcal{T} $ or $ Cn\left( \mathcal{T} \right) $;
and if one collapses every element in $ R_3 $ of $ \text{TRUTH} $ into $ R_2 $, it would clearly become a $ \mathcal{G}' $-normal structure.

Now, let $ \mathcal{G}'' = \left< \textit{D}_{\mathcal{G}}, R_i^j + \text{TRUTH} , \mathcal{P}' \right> $ be a partial pragmatic structure that includes the (partial) relation $ \text{TRUTH} = \left< R_1 , R_2 , R_3 \right> $ both in the structure and the set of primary sentences such that 
\begin{equation}
	\mathcal{P}' = \mathcal{P} \cup \left\{ \text{`TRUTH}\left( x \right) \text{'} \, \middle\vert \, x \in R_1  \right\}
	\text{ .}
\end{equation}
Therefore, from the undefinability theorem~\cite{bib11}, one will have it that there is no total expansion of $ \mathcal{G}'' $ that is a (total) model for $ \mathcal{P}' $, although we had by hypothesis that the theory $ \mathcal{T} $ was satisfiable in the same language of $ \mathcal{P}' $.
In other words, the third component $ R_3 $ in the predicate $ \text{TRUTH} = \left< R_1 , R_2 , R_3 \right> $ needs to be non-empty, whether $ R_1 $ is well-defined or not;
additionally, in case $ R_1 $ is well-defined (for example when one leaves in $ R_3 $ the problematic Gödel numbers), there will be no $ \mathcal{G}'' $-normal structure, although there was a $ \mathcal{G}' $-normal structure.
Thus, $ \left\{ \text{`TRUTH}\left( x \right) \text{'} \, \middle\vert \, x \in R_1  \right\} $ is a scheme of formulas that is formally nonsensical to add to $\mathcal{P} $ of $ \mathcal{G}' $.
\end{example}

Example~\ref{exampleFormalnonsense} highlights the importance and contribution brought about by the concept of formal nonsense for SPSs, concept which allows one to make model-theoretic assertions at the edge of what should be forbidden in model theory itself.
Although it would make sense to talk about an extra (partial) relation for actual truth relationships of a theory---because one had already assumed that there was a model for the theory in the first place---, in the moment one tries to ``mirror'' or ``force'' the insertion of this predicate into the syntactic realm (instantiated in $ \mathcal{P} $) it would lead into \emph{nonsense}.
Particularly, it is a formal nonsense because the SPS does not admit any total-structure expansion even though there already was a (total) model for the theory by construction.
At first glance, the openness and partialness in quasi-truth---unlike in classical model theory---might give the impression that such process would be possible.
However, the formal nonsense as in Definition~\ref{def:36} shows that even in scenarios in which there is a lack of complete formal knowledge, one should still be able to assert at least about the possibility of truth or falsehood after the arbitrary or subjective addition of a novel sentence in $ \mathcal{P} $.

We argue that the Example~\ref{exampleFormalnonsense} not only illustrates but also summarizes the key idea of quasi-truth in a single unary relation that represents one of the cases in which any syntactical attempt to ``grasp'' (in the theory itself) the structure along with this unary relation ultimately would prevent the existence of total expansions of such a structure.
In turn, this very unary partial relation corresponds to our notion of ``truth'' (in first order logic). 

Thus, in the context of pragmatic truth, we also argue that a sentence being formally nonsensical is a possibility that should be necessarily considered in addition to a sentence being, true, false, pragmatically true, pragmatically false.
In the next section, we will show that this concept can be extended to also encompass sentences according to a strictly pragmatic criterion.

\section{Pragmatic Nonsense}\label{sec4}

Consider, for instance, the simple pragmatic structure $\mathcal{G}$ of our Example \ref{ex:2}, but with a different set of primary sentences $\mathcal{P}_{\mathcal{G}}$ = \{‘$\alpha$ → $\beta$’\}, ‘$\alpha$ → $\beta$’ meaning in the metalanguage ‘If Joseph has at least 99,99\% of chance of being the biological father of Peter, then Joseph is the biological father of Peter’. Let us call such simple pragmatic structure $\mathcal{G}$’. $\mathcal{G}$’ may be expanded to two possible total structures: $\mathcal{H}$’, in which the partial relation \textit{$R_2^2$} is expanded to the total relation \textit{$R_2^2$}’ so that the ordered pair ⟨Joseph, Peter⟩ is in the relation; and $\mathcal{H}$’’, in which the relation is expanded so that the ordered pair ⟨Joseph, Peter⟩ is not in the relation. Both total structures are consistent with $\mathcal{P}_{\mathcal{G}}$ and may, therefore, be considered $\mathcal{A}$-normal structures.

Consider now that the DNA paternity test was made and its result was positive, that is, that the sentence ‘$\alpha$’ or ‘Joseph has at least 99,99\% of chance of being the biological father of Peter’ could have been included in $\mathcal{P}_{\mathcal{G}}$, but it was not included. Thus, one may say that the sentence ‘\textit{$P_2^2$} (\textit{Joseph}, \textit{Peter})’ or ‘Joseph is the biological father of Peter’ is pragmatically true in the simple pragmatic structure $\mathcal{G}$’, since it is true in the $\mathcal{A}$-normal structure $\mathcal{H}$’, and that the sentence ‘¬\textit{$P_2^2$} (\textit{Joseph}, \textit{Peter})’ or ‘Joseph is not the biological father of Peter’ is also pragmatically true in the simple pragmatic structure $\mathcal{G}$’, since it is true in the $\mathcal{A}$-normal structure $\mathcal{H}$’’. That is, once disregarded the observational datum constituted by the result of the DNA test, which is available, relevant and assumedly true, both hypotheses are formally admissible; nonetheless, from the pragmatic point of view, only one hypothesis is admissible: the one which takes into account the result of the DNA test.

Consider then the simple pragmatic structure $\mathcal{G}$ of our Example \ref{ex:2}, but with a different set of primary sentences $\mathcal{P}_{\mathcal{G}}$ = \{‘¬$\alpha$’\}, ‘¬$\alpha$’ meaning in the metalanguage ‘Joseph is not the biological father of Peter’, which is the verdict of a fortune teller. In this scenario, the paternity test was not made. Let us call such simple pragmatic structure $\mathcal{G}$’’. As is the case with $\mathcal{G}$, only one of the possible total expansions of $\mathcal{G}$’’ may be considered an $\mathcal{A}$-normal structure: $\mathcal{H}$’’’, in which the partial relation \textit{$R_2^2$} is expanded to the total relation \textit{$R_2^2$}’ so that the ordered pair ⟨Joseph, Peter⟩ is not in the relation. Thus, one may say that the sentence ‘¬\textit{$P_2^2$} (\textit{Joseph}, \textit{Peter})’ or ‘Joseph is not the biological father of Peter’ is pragmatically true in the simple pragmatic structure $\mathcal{G}$’’, since it is true in the $\mathcal{A}$-normal structure $\mathcal{H}$’’’. That is, once regarded such “testimonial” datum, which should, pragmatically, be considered false, but was assumed to be true, only one hypothesis is formally admissible; but, from a pragmatic point of view, the admissibility of the contradictory hypothesis, ‘\textit{$P_2^2$} (\textit{Joseph}, \textit{Peter})’, should not be excluded, for the verdict of a fortune teller does not constitute scientific evidence.

Therefore, it seems necessary that, considering the aforementioned role of the set of primary sentences $\mathcal{P}$, depending on the domain of investigation represented by a partial structure, certain restrictions or impositions are made to the inclusion of sentences in $\mathcal{P}$. That is, $\mathcal{P}$ should include the observational data and/or the theoretical results, whether true or pragmatically true, that are considered relevant and are available; and $\mathcal{P}$ should not include sentences arbitrarily assumed as true if those are inconsistent with the observational data or the theoretical results, whether true or quasi-true, that are considered relevant and are available.

Of course, a central notion here is that of \emph{relevance}: what is considered relevant or not depends not only on the subject of the investigation, but also on its context and objectives. 
Using the terminology introduced by Kuhn in \textit{The Structure of Scientific Revolutions} \cite{bib13}, the set of primary sentences $\mathcal{P}$ is not the same for an investigation carried out in a context of normal science and for a research that aims at investigating a completely new hypothesis that, if confirmed, could question the current paradigm.

For instance, a research aimed at making observations that “[...] can be compared directly with predictions from the paradigm theory” \cite[p.~26]{bib13} would consider relevant the theoretical results in question. One example would be the first direct observation of gravitational waves (see \cite{bib14}), which proved predictions made by Einstein one hundred years before the observation, within the context of the theory of general relativity.

On the other hand, a researcher may consider relevant the observational data of an anomalous phenomenon, like when Röntgen “[…] interrupted a normal investigation of cathode rays because he had noticed that a barium platino-cyanide screen at some distance from his shielded apparatus glowed when the discharge was in process” \cite[p.~57]{bib13}; that observation later led, serendipitously, to the discovery of X-rays. 

If the set of primary sentences $\mathcal{P}$ added to a partial structure is built according to those general rules, the sentences about such simple pragmatic structure may be considered pragmatically true (or pragmatically false); otherwise, it seems to us that those sentences should not be considered neither pragmatically true, nor pragmatically false, but nonsensical, since the tools used to determine their pragmatic truth or falsehood, in view of the partial character of the structure and its relations, would not constitute a reliable foundation for pragmatic truth.

Such restrictions would not allow the construction of the simple pragmatic structures $\mathcal{G}$’ and $\mathcal{G}$’’ of our examples (and, consequently, the construction of the $\mathcal{A}$-normal structures $\mathcal{H}$’, $\mathcal{H}$’’ e $\mathcal{H}$’’’), or, once they were built, would not allow one to consider the sentences about $\mathcal{G}$’ and $\mathcal{G}$’’ pragmatically true.

Thus, before defining the notion of pragmatic nonsense, we must first define when a sentence must be included in set of primary sentences $\mathcal{P}$ (and consequently when a sentence must be excluded from $\mathcal{P}$), and when the set of primary sentences $\mathcal{P}$ is well-built (and consequently when $\mathcal{P}$ is not well-built).

\begin{definition}[\textbf{Necessary inclusion and necessary exclusion of sentences to a set of primary sentences $\mathcal{P}$}]
\label{def:41}

Let $ L $ be a formal language interpreted according to an interpretation $ \mathbf{I} $ on a partial structure $ E $ with a domain $ D $, and on a simple pragmatic structure $ \mathcal{A} $ built upon $ E $. Let $\mathcal{P} $ be the set of primary sentences of $\mathcal{A} $.
\begin{itemize}
\item A sentence `$ \alpha $' of $ L $ must be \textit{necessarily included} in $ \mathcal{P} $ if, and only if, it is a true sentence relevant to $ \mathcal{A} $ or a pragmatically true sentence relevant to $\mathcal{A} $.

\item A sentence `$ \alpha $' must be \emph{necessarily excluded} from $ \mathcal{P} $ if, and only if: either `$ \alpha $' is inconsistent (that is, contradictory) with a necessarily included sentence of $ \mathcal{P} $; or `$ \alpha $' is irrelevant to $ \mathcal{A} $; or `$\alpha $' is formally nonsensical (according to Definition \ref{def:36}).

\end{itemize}
\end{definition}

\begin{definition}[\textbf{Well-built and not well-built sets of primary sentences $ \mathcal{P} $}]
\label{def:42}

Let $ L $ be a formal language interpreted according to an interpretation $ \mathbf{I} $ on a partial structure $ E $ with a domain $ D $, and on a simple pragmatic structure $ \mathcal{A} $ built upon $ E $. Let $\mathcal{P} $ be the set of primary sentences of $\mathcal{A} $.
\begin{itemize}
\item A set of primary sentences $ \mathcal{P} $ is \emph{well-built} if, and only if, it includes every necessarily included sentence and does not include any necessarily excluded sentence.

\item A set of primary sentences $ \mathcal{P} $ is \emph{not well-built} if, and only if, it does not include any necessarily included sentence or it includes any necessarily excluded sentence.

\end{itemize}

\end{definition}

Hence, we define the notion of pragmatic nonsense and \emph{strictly pragmatic} truth as follows.

\begin{definition}[\textbf{Pragmatic nonsense for formal languages interpreted on simple pragmatic structures}]
\label{def:43}

Let $ L $ be a formal language interpreted according to an interpretation $ \mathbf{I} $ on a partial structure $ E $ with a domain $ D $, and on a simple pragmatic structure $\mathcal{A}$ built upon $ E $. Let $\mathcal{P}$ be the set of primary sentences of $\mathcal{A}$.

\begin{itemize}
	\item A sentence ‘$\alpha$’ of $ L $ 
	is \textit{pragmatically nonsensical} according to the interpretation $ \mathbf{I} $ on $\mathcal{A}$ if, and only if, 
	$ \mathcal{P} $ is \emph{not} well-built 
	or ‘$\alpha$’ \emph{is} formally nonsensical.
\end{itemize}

\end{definition}

That is, every sentence of a formal language interpreted on a simple pragmatic structure whose set of primary sentences $\mathcal{P}$ is not well-built, as well as any formally nonsensical sentence, are considered pragmatically nonsensical and, therefore, not pragmatically true nor false.

Finally, following from the previous studied definitions and examples, the concept of pragmatic nonsense lead to a version of pragmatic truth in which the ``pragmaticness'' is strict, i.e., it is strictly dependent on empirical or other characteristics derived from the experience, as relevance becomes a necessary condition in the following definition:

\begin{definition}[\textbf{Strictly pragmatic truth}]
\label{def:44}

Let $ L $ be a formal language interpreted according to an interpretation $ \mathbf{I} $ on a partial structure $ E $ with a domain $ D $, and on a simple pragmatic structure $\mathcal{A}$ built upon $ E $. Let $\mathcal{P}$ be the set of primary sentences of $\mathcal{A}$.

\begin{itemize}
	\item A sentence ‘$\alpha$’ of $ L $ 
	may be pragmatically true or pragmatically false according to the interpretation $ \mathbf{I} $ on $\mathcal{A}$ if, and only if, ‘$\alpha$’ is not pragmatically nonsensical.
\end{itemize}


\end{definition}

\section{Conclusion}\label{sec5}

When constructing a formal language, the logician usually does not care about the case in which the rules are not followed. For the work of the logician, the sequences of symbols that are not formulas are not relevant. Thus, defining nonsense as those expressions that are devoid of meaning may seem unnecessary from a formal point of view, though certainly not from a philosophical one.

It seems that, in that same spirit, the theory of pragmatic truth, as developed by da Costa and his collaborators, takes for granted the best scientific practices, as it expects the relevant empirical and theoretical data to be added to a simple pragmatic structure’s set of primary sentences ${\mathcal{P}}$, disregarding the misuse of the tools it provides as irrelevant. However, we believe that, even though there can be no completely formal-theoretic criteria for the inclusion or exclusion of sentences to ${\mathcal{P}}$, there may be clearer pragmatic criteria.

Thus, our definition of pragmatic nonsense for formal languages interpreted on simple pragmatic structures (that is, Definition \ref{def:43}) is an extension of our definition of formal nonsense for formal languages interpreted on simple pragmatic structures (that is, Definition \ref{def:36}); which in turn is an extension of our definition of formal nonsense for formal languages interpreted on total structures (that is, Definition \ref{def:21}). Therefore, the notion of pragmatic nonsense extends the notion of formal nonsense, which we could call a Tarskian notion of nonsense, the same way as the theory of pragmatic truth extends the Tarskian theory of truth. It also allows one to further explore the theory of pragmatic truth, through the definition of strictly pragmatic truth, making it clear that it is not the case that any sentence about any partial structure can be considered pragmatically true simply with the addition of an empty or completely arbitrary or irrelevant set of primary sentences to such structure.

\section{Acknowledgements}

Ricardo P. Cavassane acknowledges the  support from the National Council for Scientific and Technological Development (CNPq), Brazil, grant $150915/2024$-$1$.
Felipe S. Abrah\~{a}o acknowledges support from the S\~{a}o Paulo Research Foundation (FAPESP), grants $2021$/$14501$-$8$ and $2023$/$05593$-$1$.




\bibliographystyle{sn-nature}
\bibliography{references}

\end{document}